\documentstyle[12pt]{article}
\begin{document}
\begin{center}
{\LARGE    {\bf A Bernstein theorem for complete
 spacelike constant mean curvature
 hypersurfaces in Minkowski space}} 
\end{center}
\vspace{0.6 cm}

\begin{center}
 {Huai-Dong Cao \hspace{1.0 cm} Ying Shen \hspace{1.0 cm}  Shunhui
Zhu}
\end{center}
\footnotetext[3]{The first author's research was supported in part by NSF grant \#DMS-9504925. The third author's research was supported in part by NSF grant \#DMS-9404263}

\begin{abstract}
We obtain a gradient estimate for the Gauss maps from complete spacelike
constant mean curvature hypersurfaces in Minkowski space into the hyperbolic
space. As  applications, we prove a Bernstein theorem which says that if the
image of the Gauss map is bounded from  one side, then the spacelike constant
mean curvature hypersurface must be linear. This result extends the previous
theorems obtained by B. Palmer
\cite{Pa} and Y.L. Xin
\cite{Xin1} where they assume that the image of the Gauss map is bounded.  We
also proved a Bernstein theorem for spacelike complete surfaces with parallel
mean curvature vector in four-dimensional spaces.
\end{abstract}
\vspace{1 cm}

\section{Introduction}

The classical Bernstein theorem \cite{B}
states that the only entire solution of the
minimal hypersurface equation in Euclidean space $I\!\!R^{n+1}$
\begin{eqnarray}
  \sum_{i=1}^{n}\frac {\partial}{\partial x_i}( \frac { \frac
{\partial u}{\partial x_i}}{\sqrt{1+ |\nabla u|^2}})=0 
\label{1}
\end{eqnarray}
is  linear for $n=2$. The higher dimensional version of the Bernstein
theorem was proved through the efferts of Federer, Fleming \cite{F}, De
Giorgi \cite{DG}, Almgren \cite{A}, Simons \cite{SJ}. They showed 
that the only entire solution of Eq.(\ref{1}) is linear for $n\leq 7$
while
Bombieri-De Giorgi-Giusti \cite{BDG} showed that there are nontrivial
solutions for $n>7$. In 1965, Chern
\cite{Ch} proved that the only entire solution of constant mean curvature 
hypersurface equation in $I\!\!R^{m+1}$
\begin{eqnarray*}
  \sum_{i=1}^{m}\frac {\partial}{\partial x_i}( \frac { \frac
{\partial u}{\partial x_i}}{\sqrt{1+ | \nabla u|^2}})=mH 
\end{eqnarray*}
is also linear for $2\leq m\leq 7$. Then he proposed the
generalized Bernstein problem by   
 studying  the distribution of normal vectors to a complete constant mean
curvature hypersurface. A nice result in this direction is a
theorem due to Hoffman,
Osserman, and Schoen \cite{H-O-S} which states that if the normal vectors
to a complete constant mean curvature surface in $I\!\! R^3$ lie in an
open hemisphere, then it is a plane.  

In 1968, Calabi \cite{Ca} proposed to study a similar problem in the
Minkowski space. He considered the  maximal spacelike hypersurface equation
\begin{eqnarray}
 \sum_{i=1}^{m}\frac {\partial}{\partial x_i}( \frac { \frac
{\partial u}{\partial x_i}}{\sqrt{1 - | \nabla u|^2}})=0 
\label{2}
\end{eqnarray}
 in the
Lorentz-Minkowski space $I\!\! R^{m+1}_{1}$ with coordinate $(x_0, x_1, \cdots,
x_m)$ and metric $ds^2=-(dx_0)^2 +\sum_{i=1}^{m} (dx_i)^2$. The Lorentz
metric induces a Riemannian metric on the graph of $u$   if
$|\nabla u|<1$. In \cite{Ca}, Calabi showed that, for $2\leq m \leq 4$, 
the graph of any solution to Eq.(\ref{2})  has to be a
hyperplane.  Later, Cheng and Yau \cite{C-Y} proved that 
the same conclusion holds for $m>4$. In fact they proved the parametric
version of this theorem: The only maximal spacelike hypersurface which is
a closed subset of the Lorentz-Minkowski space (with respect to the
Euclidean topology) is a linear hyperplane.

Note that a maximal spacklike hypersurface is a spacelike hypersurface
with zero constant mean curvature. In \cite{Tr}, Treibergs showed that the
constant mean curvature spacelike hypersurface equation
\begin{eqnarray*}
 \sum_{i=1}^{m}\frac {\partial}{\partial x_i}( \frac { \frac
{\partial u}{\partial x_i}}{\sqrt{1 - | \nabla u|^2}})=mH, 
\end{eqnarray*}
where $|\nabla u|<1$ and $H>0$, has many solutions besides hyperboloids.
One important fact in the study of constant mean curvature spacelike hypersurfaces of Minkowski space is that the Gauss map of such a hypersurface $M$  is a harmonic map to
hyperbolic space. Let us first recall how the Gauss map is defined. If we let $\nu$ be the timelike unit normal
vector field to $M$ in $I\!\!
R^{m+1}_{1}$. For any point $p\in M$, $|\nu (p)|^2=-1$. By parallel
translation to the origin in $I\!\! R^{m+1}_{1}$ we can regard $\nu (p)$ as a
point in $m$-dimensional hyperbolic space 
$I\!\!H^m(-1)$ 
which is
canonically embedded in $I\!\! R^{m+1}_{1}$. The construction gives us a
Gauss map  $f : M \rightarrow I\!\!H^m(-1)$. 
In \cite{Cho-Tr}, Choi and Treibergs  did an extensive study of the
interplay of properties of such a hypersurface $M$ and behavior of its Gauss map. As a result, they constructed various harmonic maps from $I\!\!R^2$ 
to $I\!\!H^2$ and $I\!\!H^2$ to $I\!\!H^2$. 

Naturally, one can ask under what conditions on the Gauss map
 a complete spacelike hypersurfaces with constant mean curvature has the
Bernstein property.  
In 1990, Palmer \cite{Pa} proved that
there exists a constant $\tau=\tau (n, H)>0$ such that if the Gauss image of
the hypersurface
$f(M)$ is contained in a geodesic ball of radius $\tau_1 < \tau$ in 
$I\!\!H^m(-1)$, then
$M$ is a hyperplane. Later, Xin \cite{Xin1} extended Palmer's theorem by
proving that
$M$ must be a hyperplane if the image under the Gauss map  is bounded.  

Our main result in this paper  is to find some natural geometric condition
and to  further extend Palmer and Xin's  Bernstein theorem by only assuming
that the image under the Gauss map is bounded from one side.  The main theorem
in Choi and Treibergs' paper
\cite{Cho-Tr} states that  given  an arbitrary closed set in the ideal
boundary at infinity of hyperbolic space, there exist a complete entire
spacelike  constant mean curvature hypersurface whose Gauss map is a
diffeomorphism onto the interior of the hyperbolic space convex hull of the
set.  Therefore, if the closure of the image under  the Gauss map intersects
the ideal boundary of hyperbolic space at more than one point,
Choi-Treibergs' theorem implies that the Bernstein theorem is not true. 
Therefore, it is natural to conjecture that the Bernstein theorem is true if
the closure of the image under  the Gauss map  intersects
the ideal boundary of hyperbolic space at exactly one point. The most
natural geometric candidate whose closure intersects the ideal boundary at
exactly one point is the so called "horoball ball" which is defined in the
next section.  If we use the upper half-plane model for
$I\!\!H^m(-1)$ so that $I\!\!H^m(-1)$ has coordinate $(\gamma_1, \cdots,
\gamma_m)$ with $\gamma_m > 0$ and the metric 
\[ds^2 =\frac { d\gamma_{1}^2 +\cdots +d\gamma_{m}^2}{ \gamma_{m}^2},\]
then the set $\{(\gamma_1, \cdots, \gamma_m) \in I\!\!H^m(-1): \gamma_m >
c >0\}$ is a horoball for any positive constant $c$. If we use the
Poinc$\acute{a}$re disk model $D$ for the hyperbolic space, then a horoball is
a smaller disk that is contained $D$ and intersects $D$ at exactly one boundary
point of $D$.  

Now,
our main theorem can be stated as follows:
\\

\noindent{\bf Main Theorem.} {\em  Let
$M$ be a complete spacelike constant mean curvature hypersurface in
Lorentz-Minkowski space
$I\!\! R^{m+1}_{1}$. Let $g$ be a function on $M$
such that $ e^{-g}$ is superharmonic. If the image of $M$ under the Gauss map $\gamma : M \rightarrow
I\!\!H^m(-1)$ lies in 
$\{ (\gamma_1, \cdots, \gamma_m):\gamma_m > e^g \}$ 
then
$M$ is necessarily a hyperplane. In particular, this is the case if the Gauss image is contained in the half-space $\{\gamma_m>
c\}$ for some constant $c>0$. }\\

Our proof of the Main Theorem is in similar spirit to that of Cheng-Yau
 \cite{C-Y} and
Xin \cite{Xin1}. The key step in the proof is  to  obtain a  gradient estimate 
for the length of the second fundamental
form of the hypersurface. This is done in section 3. The new gradient
 estimate (Theorem 3.1), which extends an earlier result of Shen \cite{Sh}, 
can be considered as the main analytic result of this paper.  It is presented
with some general curvature assumptions on both doamin and  target manifolds
because the result is of interest in its own right for the study of harmonic
maps.\\
\vspace{0.3 cm} 

{\em Acknowledgement}. The second author would like to thank Professor
Carolyn Gordon and Professor Richard Schoen for their interest in the result.

\section{Preliminaries}

In this section we shall  recall some basic facts and formulas in
Lorentzian geometry. 

Let $L$ be an $(m+1)$-dimensional Lorentzian manifold with Lorentzian metric
$\bar {g}$ of signature $(-, +, \cdots, +)$.  Let  $\{e_0, e_1, \cdots, e_m\}$
be a local Lorentzian orthonormal moving frame 
in $L$. Let $\omega_0, \omega_1, \cdots, \omega_m$ be its dual frame  so
that $\bar{g}=-\omega_0^2 + \sum_{i=1}^{m}\omega_{i}^2$. We adopt the
convention that Latin indices run from $1$ to $m$ and Greek indices run from
$0$ to $m$.  The Lorentzian connection forms $\omega_{\alpha \beta}$ of $L$
are uniquely determined by the structure equations:
\begin{eqnarray}
d\omega_0=\omega_{0i} \wedge \omega_i,  \\
d\omega_i=-\omega_{i0}\wedge \omega_0 + \omega_{ij}\wedge \omega_j,\nonumber\\
\omega_{\alpha \beta}=-\omega_{\beta \alpha}. \nonumber 
\label{st}
\end{eqnarray}

The Covariant derivatives are defined by the following equations
\begin{eqnarray}
De_0=\omega_{0i}e_i,\nonumber\\
De_i=\omega_{ij}e_j-\omega_{i0}e_0.
\label{3}
\end{eqnarray}

The curvature forms $\bar{\Omega}$ of $L$ are given by 
\begin{eqnarray*}
\bar{\Omega}_{0i}=d\omega_{0i}-\omega_{0k}\wedge \omega_{ki},\nonumber\\
\bar{\Omega}_{ij}=d\omega_{ij}+\omega_{i0}\wedge
\omega_{0j}-\omega_{ik}\wedge \omega_{kj},\\
\bar{\Omega}_{\alpha \beta}=-\frac {1}{2} \bar{R}_{\alpha \beta \gamma
\delta}\omega_{\gamma}\wedge \omega_{\delta},\nonumber
\end{eqnarray*}
where $\bar{R}_{\alpha \beta \gamma \delta}$ are components of the curvature
tensor $\bar{R}$ of $L$. 

Let $M$ be a spacelike hypersurface in a Lorentzian manifold $L$. We can
choose a local Lorentzian orthonormal  moving frame $\{e_0, e_1, \cdots, e_m\}$
in $L$ such that, restricted to $M$, the vectors $e_1, \cdots, e_m$ are
tangential to $M$ and $e_0$ is normal to $M$.  So the form $\omega_0$ is a zero
form on the hypersurface
$M$ and the induced Riemannian metric $g$ of $M$ is
\begin{eqnarray*}
g=\sum_{i} \omega_{i}^2.
\end{eqnarray*}

Hence the induced structure equations of $M$ are
\begin{eqnarray*}
d\omega_i=\omega_{ik}\wedge \omega_{k},\\
d\omega_{ij}=\omega_{ik}\wedge \omega_{kj}-\omega_{i0}\wedge \omega_{0j}
+\bar{\Omega}_{ij},\\
\omega_{ij}=-\omega_{ji},\\
\Omega_{ij}=d\omega_{ij}-\omega_{ik}\wedge \omega_{kj}=-\frac
{1}{2}R_{ijkl}\omega_k\wedge \omega_l,\\
\end{eqnarray*}
where $\Omega_{ij}$ and $R_{ijkl}$ are the curvature forms and the curvature
tensors of $M$, respectively.

Since $\omega_0=0$ on $M$, equation (\ref{st}) shows that 
\begin{eqnarray*} 
\omega_{0i}\wedge \omega_i=0.
\end{eqnarray*}

The Cartan's lemma then implies that
\begin{eqnarray*}
\omega_{i0}=h_{ij}\omega_j,
\end{eqnarray*}
where $h=h_{ij}\omega_i \otimes \omega_j$ is called the second fundamental form
of $M$. We then obtain the Gauss equation
\begin{eqnarray*}
R_{ijkl}=\bar{R}_{ijkl}-(h_{ik}h_{jl}-h_{il}h_{jk}).
\end{eqnarray*}

By taking the trace, we get the equation for the Ricci tensors 
\begin{eqnarray*}
R_{ij}=\bar{R}_{ij}+h_{ik}h_{kj}-mHh_{ij},
\end{eqnarray*}
where $H=\frac {1}{m} \sum_{i}h_{ii}$ is the mean curvature of $M$ in $L$. If
$L$ has Ricci curvature bounded from below by a constant $k_L$, then the lower
bound of the Ricci curvature of $M$ is given by
\begin{eqnarray}
 Ric_M\geq k_L -\frac {1}{4} m^2 H^2.
 \label{4}
\end{eqnarray}

We also like to introduce the notion Busemann function and its level
set in order to set up our results.  Suppose that
$N$ is a complete Riemannian manifold of nonpositive curvature and let $c:
I\!\!R
\rightarrow N$ be a unit speed geodesic. Then we call the union of balls $ B_c
=
\bigcup_{t>0}B_t(c(t))$ the horoball with center at infinity $c(\infty)$. For
$x\in N$, the function $t
\rightarrow t-d(x, c(t))$ is bounded from above and monotonically increasing.
So we can define the Busemann function by
\[ B(x)=\lim_{t\rightarrow +\infty}(t-d(x, c(t))). \]
It is known (see \cite {He-Ho} and \cite {B-Gr-S}) that the Busemann function
$B$ is a concave
$C^2$-function with
$|\nabla B|=1$.

We end this section by providing a proof of the well-known
result  that  the Gauss map of a constant mean curvature spacelike
hypersurface is a harmonic map into hyperbolic space (see
\cite{Mi},
\cite{Ish}). The relation between the second fundamental form
and the derivative of the Gauss map will become clear in the computation.

The Gauss map is a map defined by
\[ f=e_0: M\rightarrow H^m.\]
We can think of $\{e_i\}$ as an orthonormal
basis for $H^m$ at the point $e_0$. Let
$\{\theta _k\}$ be the dual basis for $H^m$ and $\{\theta _{kl}\}$ be its
connection form. Then we have
\[ f^* \theta_k=\omega_{0k}=h_{ik}\omega_i.\]

Hence we see that 
\begin{eqnarray}
f^{k}_{i}=h_{ik}.
\label {5}
\end{eqnarray}

This implies that 
\begin{eqnarray}
 |h|^2=|\nabla f|^2. 
\label{6}
\end{eqnarray}

By the second structure equation, we have 
\[ d\omega_{k0}=-\omega_{kl}\wedge \omega_{l0}.\]

So we know that 
\[ f^*\theta_{ij}=\omega_{ij}.\]

Differentiating (\ref{5}), we get 
\begin{eqnarray*}
f^{k}_{ij}\omega_j+f^{kj}\omega_{ji}-f^{li}f^*\theta_{kl}=df^{k}_{i}\\
=dh_{ik}=h_{kil}\omega_j+h_{kj}\omega_{ji}+h_{li}\omega_{lk}.
\end{eqnarray*}

Hence we obtain that 
\[ f^{k}_{ij}=h_{kij}. \]

Taking the trace, we get
\[ \sum_i f^{k}_{ii}=\sum_i h_{kii}=nH_k, \]
by the symmetry of the second fundamental form $h$. So $f$ is 
harmonic if and only if the mean curvature $H$ is constant.

\section{The Proof of Main Theorem}

In this section, we prove the Main Theorem  stated in the introduction. 
The main analytic tool involved is the following gradient estimate for 
harmonic maps. \\

\noindent{\bf Theorem 3.1.}  {\em Suppose that $M$ is a complete Riemannian
manifold with Ricci curvature $Ric_M \geq -k^2$ for some $k\geq 0$ and $N$ is
a simply-connected complete Riemannian manifold with sectional
curvature $K_N \leq -1$.  If
$f: M \rightarrow N$ is a harmonic map and $B$ is a Busemann function on $N$
such that $B\circ f \geq g$ for some function $g$ satisfying $\triangle
e^{-g}\leq 0$ (*), then we have}
\begin{eqnarray}
sup_{B_{\frac {a}{2}}(x_0)}\frac {|\nabla f|}{B\circ f - g} \leq c(m,n) \frac
{1+ka}{a}
\label {gradient}
\end{eqnarray}
{\em where $a>0$ is an arbitrary constant}.\\

From Theorem 3.1,  one can also obtain the following estimate which can be
regarded as an  extension of Yau's  Liouville Theorem for harmonic
functions.\\

\noindent{\bf Theorem 3.2.} {\em Let $f$ be a harmonic function on a complete
Riemannian manifold $M$  with  Ricci curvature
$Ric_M\geq -k^2$ for some constant $k\geq 0$. If there exists a
 function $g$ such that $e^{-g}\leq 0$ and  $f\geq g$, then} 
\begin{eqnarray*}
sup_{B_{\frac {a}{2}}(x_0)}\frac {|\nabla f|}{f- g} \leq c(m) \frac
{1+ka}{a}.
\end{eqnarray*}
{\em for all positive constant $a>0$. In particular, if the Ricci curvature
of $M$ is nonnegative, then the harmonic function has to be constant.}\\

\noindent{\bf Proof of Theroem 3.1.} Let  $M$ be a complete Riemannian
manifold with Ricci curvature $Ric_M \geq -k^2$ for some $k\geq 0$ and $N$ be
a simply-connected complete Riemannian manifold with sectional
curvature $K_N \leq -1$.  If
$f: M \rightarrow N$ is a harmonic map and $B$ is a Busemann function on $N$
such that $B\circ f \geq g$ for some function $g$ satisfying $\triangle
e^{-g}\leq 0$. Let 
\[\phi=\frac {|\nabla f|}{B\circ f -g},\]
then
we have 
\begin{eqnarray*}
\nabla \phi = \frac {\nabla |\nabla f|}{B\circ f-g}-\frac { |\nabla
f|(\nabla B\circ f- \nabla g)}{ (B\circ f -g)^2}
\end{eqnarray*}
and

\begin{eqnarray}
\triangle \phi = \frac {\triangle |\nabla f|}{B\circ f -g}-2\frac {\nabla
(B\circ f- g)\cdot \nabla \phi}{B\circ f -g} -\frac { \phi \triangle B\circ f -
\phi \triangle g}{B\circ f -g}.
\label{11}
\end{eqnarray}

For harmonic map $f$, we have the following Bochner formula
\begin{eqnarray}
\frac {1}{2}\triangle |\nabla f|^2 \geq |f^{k}_{ij}|^2-k^2|\nabla f|^2
\label{12}
\end{eqnarray}

It is known that (see \cite{Sc-Yau})
\begin{eqnarray}
|f^{k}_{ij}|^2 \geq (1+ \frac {1}{2mn})|\nabla |\nabla f||^2,
\label{13}
\end{eqnarray}
where $m=dim(M)$ and $n=dim (N)$. Combining (\ref{12}) and (\ref{13}), we
obtain
\begin{eqnarray}
|\nabla f|\triangle |\nabla f| \geq \epsilon |\nabla |\nabla f||^2 -
k^2|\nabla f|^2,
\label{14}
\end{eqnarray}
where $\epsilon = \frac {1}{2mn}$.

Inequalities (\ref{14}) and (\ref{11}) imply
\begin{eqnarray}
\triangle \phi \geq \epsilon \frac {|\nabla |\nabla f||^2}{|\nabla
f|(B\circ f-g)} -k^2\phi \nonumber\\
 -2\frac {(\nabla B\circ f-\nabla g) \cdot \nabla
\phi}{B\circ f -g} -\frac {\phi \triangle ( B\circ f -g)}{B\circ f -g}.
\label{15}
\end{eqnarray}

We know that 
\begin{eqnarray}
& &-2 \frac {(\nabla B\circ f-\nabla g) \cdot \nabla \phi}{B\circ f
-g} \nonumber\\
&=& -(2-2\epsilon) \frac {(\nabla B\circ f - \nabla g)\cdot
\nabla \phi}{B\circ f -g}  -2\epsilon \frac {(\nabla B\circ f
-\nabla g)\cdot \nabla \phi)}{B\circ f- g} \nonumber
\\ 
&=&-(2-2\epsilon) \frac {(\nabla B\circ f - \nabla g)\cdot
\nabla \phi}{B\circ f -g} -2\epsilon \frac { (\nabla B\circ f -\nabla g)\cdot
\nabla |\nabla f|}{ (B\circ f-g)^2} \nonumber\\
&+&2\epsilon \frac {|\nabla B\circ f -\nabla g|^2|\nabla f|}{(B\circ f -g)^3 } \nonumber\\
&\geq &-(2-2\epsilon) \frac {(\nabla B\circ f - \nabla
g)\cdot
\nabla
\phi}{B\circ f -g} 
 -2\epsilon \frac {|\nabla B\circ f-\nabla g||\nabla |\nabla f||}{ (B\circ f
-g)^2} \nonumber \\
&+& 2\epsilon \frac { |\nabla B\circ f -\nabla g|^2|\nabla
f|}{(B\circ f -g)^3} ,
\label{16}
\end{eqnarray}
and
\begin{eqnarray}
& &-2\epsilon \frac {|\nabla B\circ f-\nabla g||\nabla |\nabla f||}{
(B\circ f -g)^2} \nonumber \\
&\geq &-\epsilon \frac {|\nabla |\nabla
f||^2}{(B\circ f-g)|\nabla f|} 
+ \epsilon \frac { |\nabla B\circ f -\nabla g|^2 |\nabla
f|}{ (B\circ f -g)^3 }.
\label{17}
\end{eqnarray}

Inequalities (\ref{15}), (\ref{16}) and (\ref{17}) then imply that
\begin{eqnarray}
\triangle \phi \geq -(2-2\epsilon)\frac {\nabla (B\circ f -g)\cdot \nabla
\phi}{B\circ f -g} + \nonumber \\
\epsilon \frac { |\nabla f||\nabla B\circ f
-\nabla g|^2}{(B\circ f -g)^3} -k^2 \phi  - \frac {\phi \triangle B\circ f-\phi
\triangle g}{ B\circ f -g}.
\label{18}
\end{eqnarray}

For any fixed point $x_0 \in M$, we can define a function $F$ on the ball
$B_a(x_0)$ by
\begin{eqnarray*}
F(x)=(a^2 -r^2)\phi (x),
\end{eqnarray*}
where $r(x)=dist(x, x_0)$ is the distance function on $M$. Without lossing of
generality, we may assume $r\in C^2(M)$.  If $\nabla f$ is not identically
zero, then $F$ must attain its maximum at some interior point $x^*$. By the
mzximum principle, we have 
\begin{eqnarray}
\nabla F( x^*) =0,
\label{19}
\end{eqnarray}
\begin{eqnarray}
\triangle F(x^*) \leq 0.
\label{20}
\end{eqnarray}
By   (\ref{19}) and (\ref{20}), we have at $x^*$

\[ \frac {\nabla r^2}{a^2 -r^2}=\frac {\nabla \phi}{\phi},\]
and
\[ -\frac {\triangle r^2}{a^2 -r^2}+\frac {\triangle \phi}{\phi}-\frac
{2\nabla r^2 \cdot \nabla \phi}{(a^2-r^2)\phi}\leq 0. \]
It follows that 
\[ \frac {\triangle \phi}{\phi}-\frac {\triangle r^2}{a^2-r^2}-\frac
{2|\nabla r^2|^2}{(a^2-r^2)^2}\leq 0.\]
Using the Laplacian comparison theorem (see \cite{Gr-Wu}), we obtain
\[ \triangle r^2 \leq C(1+kr).\]
Also, it is well-known that 
\[ |\nabla r|=1.\]
Using inequality  (\ref{18}), we obtain 
\begin{eqnarray}
\epsilon \frac {|\nabla B\circ f -\nabla g|^2}{(B\circ f -g)^2 } -
\frac {\triangle (B\circ f -g)}{ B\circ f -g} - \frac {4(1-\epsilon)r}{a^2
-r^2}\frac {\nabla (B\circ f -g)\cdot \nabla r}{B\circ f -g}\leq A,
\label{21}
\end{eqnarray}
where 
\[ A=k^2 + \frac {C (1+kr)(a^2 -r^2)+8r^2}{(a^2-r^2)^2}. \]
Using $ab\geq -\frac {\epsilon}{2}a^2 -\frac {2}{\epsilon}b^2$ and $|\nabla
r|=1$,   (\ref{21}) can be simplified as 
\begin{eqnarray}
\frac {\epsilon}{2} \frac {|\nabla B\circ f -\nabla g|^2}{(B\circ f -g)^2 } -
\frac {\triangle (B\circ f -g)}{ B\circ f -g} \leq A',
\label{22}
\end{eqnarray}
where 
\[ A'=k^2 + C\frac { (1+kr)(a^2 -r^2)+r^2}{(a^2-r^2)^2}. \]
It is easy to see that

\begin{eqnarray}
|\nabla B\circ f -\nabla g|^2\geq \frac {1}{2} |\nabla B\circ f|^2 -|\nabla
g|^2.
\label{23}
\end{eqnarray}

Since $B$ is a horofunction, we can always assume that $B\circ f \geq 1+ g$.
Combining (\ref{22}) and (\ref{23}), we have
\begin{eqnarray}
\frac {\epsilon}{4} \frac {|\nabla B\circ f|^2}{(B\circ f -g)^2}-\frac
{\triangle B\circ f}{B\circ f -g} + \frac {\triangle g}{B\circ f -g} - 
\frac {\epsilon}{2}\frac {|\nabla g|^2}{ (B\circ f -g)^2} -A'\leq 0.
\label{24}
\end{eqnarray} 

By Condition ($\ast$),  we see that 
\[ |\nabla g|^2 - \triangle g \leq 0.\]
So 
\begin{eqnarray}
\frac {\triangle g}{B\circ f -g} - 
\frac {\epsilon}{2}\frac {|\nabla g|^2}{ (B\circ f -g)^2}\geq 0.
\label{25}
\end{eqnarray}
From  (\ref{24}) and (\ref{25}), we get 
\begin{eqnarray}
\frac {\epsilon}{4} \frac {|\nabla B\circ f|^2}{(B\circ f -g)^2}-\frac
{\triangle B\circ f}{B\circ f -g} \leq A'.
\label{26}
\end{eqnarray}
If we choose a local polar coordinate chart around $p=f(x)$ such that 
$e_1=\nabla B$, the radial direction. Then we have
\begin{eqnarray}
 |\nabla B\circ
f|^2=\sum_{i=1}^{m}|f^{\alpha}_{i}B_{\alpha}(f)|^2\geq 
\sum_{i=1}^{m}|f^{1}_{i}(x)|^2.
\label{A}
\end{eqnarray}
Applying the Hessian comparison theorem, we have
\[ B_{\alpha \beta}=-\lim_{t\rightarrow \infty}d_{\alpha \beta}(x, c(t)\leq
-\lim_{t\rightarrow \infty}\coth (d(x, c(t))H_{\alpha \beta}(p)\leq -H_{\alpha
\beta}(p), \]
where $H_{\alpha \beta}=\gamma_{\alpha \beta}-dB_{\alpha}dB_{\beta}$ and
$\gamma_{\alpha \beta}$ is the metric on $N$. 

Since $f$ is harmonic, 
\begin{eqnarray}
\triangle B\circ
f(x)=g^{ij}f^{\alpha}_{i}f^{\beta}_{j}B_{\alpha
\beta}
\leq -\sum_{\alpha =2}^{n}\sum_{i=1}^{m}|f^{\alpha}_{i}|^2.
\label{B}
\end{eqnarray}

Combining (\ref{A}), (\ref{B}) and using the assumption that   $B\circ f -g\geq 1$, we obtain the  the following key inequality
\begin{eqnarray}
\frac {\epsilon}{4} \frac {|\nabla B\circ f|^2}{(B\circ f -g)^2}-\frac
{\triangle B\circ f}{B\circ f -g} \geq \frac {\epsilon}{4} \frac {|\nabla
f|^2}{(B\circ f -g)^2}.
\label {27}
\end{eqnarray}

Inequalities (\ref{26}) and (\ref{27}) then imply that 
\begin{eqnarray*}
\frac {\epsilon}{4} \frac {|\nabla
f|^2}{(B\circ f -g)^2} \leq A'.
\end{eqnarray*}
which is equivalent to 
\begin{eqnarray}
\epsilon F^2 \leq A'.
\label{30}
\end{eqnarray}

Now the  gradient estimate
(\ref{gradient}) follows easily from (\ref{30}). This completes the proof of Theorem 3.1.\\

Using Theorem 3.1, we now present the proof of the Main Theorem.\\

\noindent{\bf Proof of Main Theorem.}  We first use the generalized  
maximum principle due to Omori \cite{Om} and Yau
\cite{Yau} to prove that the constant mean curvature is actually zero. 

Let $f$ be the Gauss map from $M$ to the hyperbolic space $I\!\! H^m$ and $B$
be
 a Busemann function on $I\!\! H^m$. Define the function $h=-(B\circ f -g)$,
where $e^{-g}$ is superharmonic on $M$. Then it is easy to see that  $h\in
C^2(M)$ \cite{He-Ho} and $h$ is bounded from above.  Since the Ricci curvature
of $M$ is bounded from below by a constant, we have the Omori-Yau \cite{Om},
\cite{Yau} maximum principle for complete manifold. Therefore, for $\forall
\epsilon >0$, and $\forall x_0 \in M$, there exists a point $x\in M$ such that 
\begin{eqnarray}
h(x)\geq h(x_0), \nonumber \\
|\nabla h|(x_0)< \epsilon,
\label{8}
\end{eqnarray}

\begin{eqnarray}
(\triangle h) < \epsilon.
\label{9}
\end{eqnarray}

From our estimates (\ref{A}) and (\ref{B}) in the proof of Theorem 3.1 and our
Condition $(\ast)$, we know that 
\begin{eqnarray}
|\nabla h|^2 + \triangle h =|\nabla B\circ f|^2 + |\nabla g|^2 -2\nabla
B\circ f \nabla g -\triangle B\circ f +\triangle g \nonumber \\
\geq (\frac {1}{2}|\nabla B\circ f|^2 -\triangle
B\circ f )+ (\triangle g- |\nabla g|^2 ) \nonumber \\
\geq \frac {1}{2}|\nabla f|^2 
=\frac {1}{2}\sum_{i, j} h^2_{ij} \geq \frac {m}{2} H^2
\label{10}
\end{eqnarray}
 
Since $H$ is a constant, inequalities (\ref{8}), (\ref{9}) and (\ref{10})
force $H\equiv 0$. Now it follows from (\ref{4}) that the Ricci curvature of $M$ is nonnegative.
Hence our gradient estimate (\ref{gradient}) in Theorem 3.1, with $k=0$, implies
\[ sup_{B_{\frac {a}{2}}(x_0)}\frac {|\nabla f|}{B\circ f - g} \leq c(m,n) 
\frac
{1}{a}.
\]
So  $|\nabla
f|^2=0$ by letting 
$a \rightarrow \infty$. Hence (\ref{6}) implies that the second 
fundamental form of $M$ is identically zero. Therefore $M$ must be a
hyperplane. This completes the proof of our Main Theorem.\\

The proof of Theorem 3.2 is exactly the same as the proof of Theorem 3.1. 
For  $N=I\!\!H^n$ and any fixed constant $c>0$, if we
use the coordinate chart
$z=(z_1, \cdots, z_n)$ with $z_n>0$ and choose 
\[\gamma=\{(0, \cdots, 0, z_n): z_n\geq c \}\]
to be a geodesic ray. Then the Busemann function for $\gamma$ is 
\[ B_c(z)=\ln (\frac {z_n}{c}).\]
For harmonic function $z=f(x)$, we can let
\[ y=\ln z. \]
Then 
\[ dz^2=\frac {dy^2}{y^2}.\]

Analytically, 
we
remark that in order to  obtain an upper bound for   $\frac {|\nabla f|}{\phi
\circ f -\tilde{g}}$, it suffices to have the followings:
\begin{eqnarray}
|\nabla \phi\circ f|^2 - \triangle \phi \circ f \geq c |\nabla f|^2,
\label{31}
\end{eqnarray}
and 
\begin{eqnarray}
|\nabla \tilde{g}|^2 - \triangle \tilde{g} \leq 0.
\label{32}
\end{eqnarray}
Inequality (\ref{27}) is just a consequence of (\ref{31}) and (\ref{32}) since
we can always make $\phi \circ f -g \geq 1$.

Clearly, our assumptions in Theorem 3.2 imply both
(\ref{31}) and (\ref{32}). Therefore, Theorem 3 is true.

\section{Applications}

In this section we are going to apply our Main Theorem to study the Bernstein
problem on spacelike 
 submanifolds with parallel mean curvature vector in semi-Euclidean
space with higher codimensions. No much is known in this direction. Our
work is inspired by Xin's original paper \cite{Xin2} which proves a nice
splitting theorem for certain spacelike surface in $I\!\!R^{4}_{2}$ . 

Let $I\!\!R^{m+n}_{n}$ be a $m+n$-dimensional semi-Euclidean space with
the semi-Riemannian metric
\[ ds^2=dx_{1}^{2} +\cdots + dx_{m}^{2} - dx_{m+1}^{2}-\cdots -
dx_{m+n}^{2}.\]
Let M be a spacelike oriented $m$-dimensional submanifold of
$I\!\!R^{m+n}_{n}$. Following \cite{Xin2},  one can define the generalized
Gauss map
$
\gamma: M
\rightarrow G_{m, n}^{m}$, where $G_{m, n}^{m}$ is a Cartan-Hardamard
manifold formed by replacing complex numbers with real numbers in the bounded
symmetric domain of first type $G_{m,n}^{m}(C)$(see \cite{Xin2}). In
the following we are going to study spacelike surfaces with parallel mean
curvature vector in $I\!\! R^{4}_{2}$, a $4$-manifold with signature $(+, +, -,
-)$. \\

\noindent{\bf Lemma 4.1.}(\cite{Hel}) {\em $G^{2}_{2,2}$ is isometric to
$I\!\!H^2(-1)\times
I\!\! H^2(-1)$. }\\

In the following we will identify
$G^{2}_{2,2}$ with
$I\!\!H^2(-1)\times
I\!\! H^2(-1)$. Then for any given map $f: M \rightarrow G^{2}_{2,2}$,  
we can write $f=(f_1, f_2)$, where $f_i=\pi_i\cdot f$ and $\pi_i: G^{2}_{2,2}
\rightarrow I\!\! H^2(-1)$ is the pojection of $G^{2}_{2,2}$ into its $i$-th
factor, $i=1, 2$.

Now, we can state a
Bernstein theorem for spacelike surface with parallel mean curvature vector in
$I\!\! R^{4}_{2}$.\\

\noindent{\bf Theorem 4.2.} {\em Let $M$ be a complete spacelike surface
with parallel mean curvature vector in $I\!\! R^{4}_{2}$ and let
$\gamma=(\gamma_1, \gamma_2): M \rightarrow G^{2}_{2,2}$ be the Gauss map. If 
 the image under the map $\gamma_1$ (or $\gamma_2$) lies in a horoball in 
   $I\!\! H^2(-1)$, then $M$ must be a plane.}\\

We choose a local Lorentzian  orthonormal frame $\{e_1, e_2, e_3, e_4\}$ along
$M$ such that
$e_1, e_2$ are tangent to $M$ with dual frame $\{\omega_1, \omega_2,
\omega_3, \omega_4\}$.  In the following, we adopt the convention that the
Latin indices run from $1$ to $2$ and the Greek indices run from $3$ to $4$.
The induced Riemannian metric of
$M$ is
\[g=\omega_{1}^2 +
\omega_{2}^{2}.\]
 So the structure equations of $M$ are 
\begin{eqnarray*}
d\omega_i=\omega_{ik}\wedge \omega_k, \nonumber\\
\omega_{ik}=-\omega_{ki}, \nonumber\\
d\omega_{ij}=\omega_{ik}\wedge \omega_{kj} - \omega_{i \alpha}\wedge
\omega_{\alpha j}, \nonumber\\
\Omega_{ij}=d\omega_{ij}-\omega_{ik}\wedge \omega_{kj}=-\frac
{1}{2}R_{ijkl}\omega_{k}\wedge \omega_{l}.\nonumber \\
\end{eqnarray*}

By Cartan's lemman, we have
\[ \omega_{i\alpha}=h^{\alpha}_{ij}\omega_{j},\]
where $h^{\alpha}_{ij}$ are components of the second fundamental form of $M$
in $I\!\! R^{4}_{2}$. The mean curvature vector of $M$
in $I\!\! R^{4}_{2}$ is
\[ \vec{H}=\frac {1}{2} g^{ij}h^{\alpha}_{ij}e_{\alpha}.\]

We then obtain from the Gauss equation that the sectional curvature of $M$ is
\[ R_{ijkl}=-(h^{\alpha}_{ik}h^{\alpha}_{jl}-
h^{\alpha}_{il}h^{\alpha}_{jk}),\]
and the Ricci curvature of $M$ is
\[ R_{jl}=g^{ik}R_{ijkl}=-(h^{\alpha}_{ii}h^{\alpha}_{jl}-
h^{\alpha}_{il}h^{\alpha}_{ji}).\]

The second structure equations on the induced normal connection of $M$ in 
$I\!\! R^{4}_{2}$ are
\begin{eqnarray*}
d\omega_{\alpha \beta}=-\omega_{\alpha \gamma}\wedge \omega_{\gamma \beta} +
\Omega_{\alpha \beta},\nonumber\\
\Omega_{\alpha \beta}=-\frac {1}{2}R_{\alpha \beta \gamma
\delta}\omega_{\gamma}\wedge \omega_{\delta},\nonumber\\
R_{\alpha \beta
ij}=-(h^{\alpha}_{ki}h^{\beta}_{kj}-h^{\alpha}_{kj}h^{\beta}_{ki}.)
\nonumber\\
\end{eqnarray*}

Define the covariant derivative of $h^{\alpha}_{ij}$ by
\[ h^{\alpha}_{ijk}\omega_{k}=dh^{\alpha}_{ij} + h^{\alpha}_{kj}\omega_{ki} +
h^{\alpha}_{ik}\omega_{kj}-h^{\beta}_{ij}\omega{\alpha \beta}.\]
It is easy to see that $h^{\alpha}_{ijk}$ is symmetric in $i, j, k$. 

We say that $M$ has parallel mean curvature vector if
\begin{equation}
 D\vec{H}=\frac {1}{2} g^{ij}
h^{\alpha}_{ijk}\omega_{k}e_{\alpha}=0
\label{34}
\end{equation}
holds everywhere in $M$.

For the Gauss map 
\[ \gamma: M \rightarrow G^{2}_{2,2},\]
so we have 
\begin{equation}
 \gamma^{\ast}\omega_{\alpha i}=h^{\alpha}_{ij}\omega_j.
\label{35}
\end{equation}

It follows from  Equations (\ref{34}) and (\ref{35}) that the Gauss map
$\gamma$ is harmonic if and only if $M$ has parallel mean curvature vector.

It is computed in \cite{Xin2} that the canonical  metric on $G^{2}_{2,2}=
I\!\!H^2(-1)\times
I\!\! H^2(-1)$
is 
\[ ds_{1}^{2}+ds_{2}^{2},\]
where 
\begin{equation}
ds_{1}^{2}= \frac {1}{2}[ (\omega_{13}+\omega_{24})^2 +(\omega_{23}
-\omega_{14})^2]
\label{36}
\end{equation}
is the induced metric of the first factor.

Equations (\ref{35}) and  (\ref{36}) then imply that
\begin{eqnarray*}
\gamma_{1}^{\ast}(\omega_{13} + \omega_{24})=(h_{11}^{3}+h^{4}_{12})\omega_1
+(h^{3}_{12}+h^{4}_{22})\omega_2,\nonumber\\
\gamma_{1}^{\ast}(\omega_{23} - \omega_{14})=(h_{11}^{4}-h^{3}_{12})\omega_1
+(h^{4}_{12}-h^{3}_{22})\omega_2.\nonumber
\end{eqnarray*}
If we write
\begin{eqnarray*}
\gamma_{1}^{\ast}(\omega_{13} + \omega_{24})=a_{11}\omega_1
+a_{12}\omega_2,\nonumber\\
\gamma_{1}^{\ast}(\omega_{23} - \omega_{14})=a_{21}\omega_1
+a_{22}\omega_2,\nonumber
\end{eqnarray*}
then we have
\begin{eqnarray*}
a_{11}=h^{3}_{11}+h^{4}_{12},\nonumber\\
a_{12}=h^{3}_{12}+h^{4}_{22},\nonumber\\
a_{21}=h^{4}_{11}-h^{3}_{12},\nonumber\\
a_{22}=h^{4}_{12}-h^{3}_{22}.\nonumber
\end{eqnarray*}

Consequently,
\begin{eqnarray*}
a_{111}=h^{3}_{111}+h^{4}_{121},\nonumber\\
a_{112}=h^{3}_{112}+h^{4}_{122},\nonumber\\
a_{121}=h^{3}_{121}+h^{4}_{221},\nonumber\\
a_{122}=h^{3}_{122}+h^{4}_{222},\nonumber\\
a_{211}=h^{4}_{111}-h^{3}_{121},\nonumber\\
a_{212}=h^{4}_{112}-h^{3}_{122},\nonumber\\
a_{221}=h^{4}_{121}-h^{3}_{221},\nonumber\\
a_{222}=h^{4}_{122}-h^{3}_{222}.\nonumber
\end{eqnarray*}

Thus, the map $\gamma_1$ is  harmonic if and only if
\begin{eqnarray}
h^{3}_{111}+h^{3}_{122}+h^{4}_{211}+h^{4}_{222}=0,\nonumber\\
h^{4}_{111}+h^{4}_{122}-h^{3}_{211}-h^{3}_{222}=0.
\label{37}
\end{eqnarray}

Similarly, the map $\gamma_2$ is  harmonic if and only if
\begin{eqnarray}
h^{3}_{111}+h^{3}_{122}-h^{4}_{211}-h^{4}_{222}=0,\nonumber\\
h^{4}_{111}+h^{4}_{122}+h^{3}_{211}+h^{3}_{222}=0.
\label{38}
\end{eqnarray}

If $M$ has parallel mean curvature vector, then 
\[ D\vec{H}=\vec{0}.\]
So $M$ has parallel mean curvature vector if and only if
\begin{eqnarray}
h^{3}_{111}+h^{3}_{122}=0,\nonumber\\
h^{3}_{211}+h^{3}_{222}=0,\nonumber\\
h^{4}_{111}+h^{4}_{122}=0, \nonumber\\
h^{4}_{211}+h^{4}_{222}=0, 
\label{39}
\end{eqnarray}

Therefore, it is easy to see from the Equations (\ref{37}), (\ref{38}) and
(\ref{39}) that a spacelike surface in $I\!\!R^{4}_{2}$ has parallel mean
curvature vector if and only if both $\gamma_1$ and $\gamma_2$ are harmonic
maps.

In the following, we assume that $M$ has parallel mean curvature vector and we
specify our local orthonormal frame
$\{ e_1, e_2, e_3, e_4
\}$ so that
\[ e_3=\frac { \vec{H}}{|\vec{H}|}.\]
Since $e_3$ is parallel in the normal bundle, we know that
\[ \omega_{34}=0,\]
and $h^{\alpha}_{ij}$ can be diagonalized simultaneously.  Therefore, we have
\begin{eqnarray*}
a_{11}=h^{3}_{11},\nonumber\\
a_{12}=h^{4}_{22},\nonumber\\
a_{21}=h^{4}_{11},\nonumber\\
a_{22}=-h^{3}_{22}.\nonumber
\end{eqnarray*}

Thus, the energy density of $\gamma_1$ is
\[  e(\gamma_1)=\frac{1}{2}\sum _{ij} a_{ij}^{2}=\frac {1}{2}[ (h^{3}_{11})^2
+(h^{3}_{22})^2+(h^{4}_{11})^2 +(h^{4}_{22})^2]. \]
Similarly, we have
\[  e(\gamma_2)=\frac {1}{2}[ (h^{3}_{11})^2
+(h^{3}_{22})^2+(h^{4}_{11})^2 +(h^{4}_{22})^2]=e(\gamma_1). \]

If the image under the harmonic map
$\gamma_1$ lies in some horoball in
$I\!\!H^2(-1)$, then we can apply our  gradient estimates in section 3 to
show that 
\[ e(\gamma_1)=0,\]
and, 
\[ e(\gamma)=e(\gamma_1)+ e(\gamma_2)=2e(\gamma_1)=0.\]
Therefore, $\gamma$ is a constant map and consequently $M$ must be a plane.
This completes our proof for Theorem 4.2.

\noindent Department of Mathematics, Texas A$\&$M University, College
Station, TX 77843\\
{\em E-mail address}: cao@math.tamu.edu\\

\noindent Department of Mathematics, Dartmouth College, Hanover, NH 03755\\
{\em E-mail address}: ying.shen@dartmouth.edu\\

\noindent Department of Mathematics, Dartmouth College, Hanover, NH 03755\\
{\em E-mail address}: shunhui.zhu@dartmouth.edu

\end{document}